\input amstex
\input amsppt.sty

\TagsOnRight

\define\Y{\Bbb Y}
\define\Z{\Bbb Z}
\define\C{\Bbb C}
\define\R{\Bbb R}

\define\be{\beta}

\define\Ga{\Gamma}
\define\de{\delta}

\define\la{\lambda}

\define\th{\theta}

\define\om{\omega}

\define\w{\tfrac{\xi}{\xi-1}}

\define\X{\frak X}

\define\ti{\tilde}
\define\wt{\widetilde}
\define\wh{\widehat}
\define\tht{\thetag}

\define\Prob{\operatorname{Prob}}

\define\const{\operatorname{const}}

\define\Conf{\operatorname{Conf}}

\define\unX{\underline X\,}
\define\unK{\underline K\,}
\define\hatK{\underline{\wh K}\,}

\define\m{{\text{Meixner}}}
\define\k{{\text{Krawtchouk}}}

\define\Me{\frak M}

\topmatter
\title Meixner polynomials and random partitions
\endtitle

\author
Alexei Borodin and Grigori Olshanski
\endauthor

\dedicatory Dedicated to our teacher A.~A.~Kirillov on the occasion of his 70th
birthday\enddedicatory

\abstract The paper deals with a 3--parameter family of probability measures on
the set of partitions, called the z--measures. The z--measures first emerged in
connection with the problem of harmonic analysis on the infinite symmetric
group. They are a special and distinguished case of Okounkov's Schur measures.
It is known that any Schur measure determines a determinantal point process on
the 1--dimensional lattice. In the particular case of z--measures, the
correlation kernel of this process, called the discrete hypergeometric kernel,
has especially nice properties. The aim of the paper is to derive the discrete
hypergeometric kernel by a new method, based on a relationship between the
z--measures and the Meixner orthogonal polynomial ensemble. In another paper
(Prob. Theory Rel. Fields 135 (2006), 84--152) we apply the same approach to a
dynamical model related to the z--measures.
\endabstract

\endtopmatter
\document

\head Introduction
\endhead

\subhead{Main definitions and motivations}\endsubhead Recall that a {\it
partition\/} is an infinite monotone sequence of nonnegative integers,
$\la=(\la_1\ge\la_2\ge\dots)$, with finitely many nonzero terms $\la_i$. There
is a natural identification of partitions with Young diagrams; for this reason,
we denote the set of all partitions by symbol $\Y$. Clearly, $\Y$ is a
countable set. To each partition $\la\in\Y$ we assign a weight depending on
three parameters $z$, $z'$, and $\xi$. Under suitable restrictions on the
parameters (for instance, if $z$ and $z'$ are complex numbers conjugate to each
other and $0<\xi<1$) all the weights are nonnegative and their sum equals 1.
Then we get a probability measure on the set $\Y$, which makes it possible to
speak about {\it random\/} partitions. The measures on $\Y$ obtained in this
way are called the {\it z--measures\/} and denoted as $M_{z,z',\xi}$ (see
section 1 for precise definitions).

Our interest in the z--measures is mainly motivated by the fact that they play
a crucial role in harmonic analysis on the infinite symmetric group, see
\cite{KOV1}, \cite{KOV2}, \cite{BO2}, \cite{Ol}. On the other hand, for special
values of parameters $z,z'$ the z--measures turn into discrete orthogonal
polynomial ensembles which in turn are related to interesting probabilistic
models: the directed percolation model \cite{Jo1}, the stochastic growth model
of \cite{GTW}, random standard tableaux of rectangular shape \cite{PR}. The
z--measures are studied in many research papers: \cite{BO2}, \cite{BO3},
\cite{BO4}, \cite{BO5}, \cite{BO6}, \cite{BOS}, \cite{Ok2} ; see also the
expository papers \cite{BO1}, \cite{Ol}. Finally, note that the z--measures are
a particular case of more general objects, the {\it Schur measures\/}
introduced by Okounkov in \cite{Ok1} and further investigated by many people.

Although the z--measures are quite interesting by themselves, the main problems
concern their limits as parameter $\xi$ approaches the critical value 1
(parameters $z,z'$ being fixed). Note that, as $\xi\to1$,  the weight of each
partition tends to 0, that is, the measure runs away to infinity. Thus, to
catch possible limits we have to embed $\Y$ in a larger space. It turns out
that there are different limit regimes, and for each regime the limit measure
lives on a suitable space of infinite point configurations (see our paper
\cite{BO5} for more details). In other words, the limit measure determines a
{\it random point process\/}. An appropriate way to describe point processes is
to use the language of correlation functions, and the first necessary step is
to interpret the initial z--measures as point processes, too.

To this end, we use a well--known interpretation of partitions as {\it Maya
diagrams\/}, which are semi--infinite point configurations on the
1--dimensional lattice. It is convenient to identify the lattice with the
subset $\Z':=\Z+\frac12\subset\R$ of (proper) half--integers. Then the Maya
diagram of a partition $\la\in\Y$ is the configuration (or simply the subset)
$\{\la_i-i+\frac12\mid i=1,2,\dots\}\subset\Z'$. Each z--measure $M_{z,z',\xi}$
thus gives rise to a {\it random\/} point configuration on $\Z'$ (or a point
process on $\Z'$), and its $n$th {\it correlation function\/} $\rho_n$
($n=1,2,\dots$) expresses the probability $\rho_n(x_1,\dots,x_n)$ that the
random configuration contains an arbitrary prescribed finite set of points
$x_1,\dots,x_n$ in $\Z'$.

It is worth noting that the correlation functions survive in various limit
regimes, which explains their efficiency.

A remarkable property of the z--measures is that, for any $n=1,2,\dots$, the
probability $\rho_n(x_1,\dots,x_n)$ can be written as the $n\times n$
determinant $\det[K(x_i,x_j)]$ where $K(x,y)$ is a function on $\Z'\times\Z'$
not depending on $n$ (it depends on parameters $z,z',\xi$ only). Random point
processes with such a property are called {\it determinantal\/}, \footnote{This
term, introduced in \cite{BO2} and then employed in Soshnikov's expository
paper \cite{S}, is now widely used. Earlier works used the term ``fermion point
processes''.} and the function $K(x,y)$ is called the {\it correlation
kernel\/}.

As was first shown in \cite{BO2}, the correlation kernel of the z--measure
$M_{z,z',\xi}$ can be explicitly written in terms of the Gauss hypergeometric
functions; for this reason we called it the {\it discrete hypergeometric
kernel\/}. Then a number of different proofs were suggested in \cite{Ok1} (see
also \cite{BOk}), \cite{Ok2}, \cite{BOS}. The goal of the present paper is to
better understand the nature of this kernel.

\subhead{The results}\endsubhead Now we are in a position to describe our main
results:
\smallskip

(1) We introduce a system of functions $\psi_a(x)=\psi_a(x;z,z',\xi)$, where
the triple $(z,z',\xi)$ is the parameter of the z--measure, $ x$ is the
argument ranging over $\Z'$, and $a\in\Z'$ is an additional parameter. For
fixed $(z,z',\xi)$ and varying $a$, the family $\{\psi_a\}$ forms an orthogonal
basis in the coordinate Hilbert space $\ell^2(\Z')$. Each function $\psi_a$ can
be expressed through the Gauss hypergeometric function.
\smallskip
(2) We exhibit a second order difference operator $D=D(z,z',\xi)$ on $\Z'$
which is diagonalized in the basis $\{\psi_a\}$. The eigenvalue of $D$
corresponding to the eigenfunction $\psi_a$ is equal to $a(1-\xi)$. (We assume
$0<\xi<1$, so that the eigenvalue $a(1-\xi)$ is positive or negative depending
on the sign of parameter $a\in\Z'$.)
\smallskip
(3) Set $\Z'_+=\{\frac12, \frac32,\frac52,\dots\}$. We prove that the discrete
hypergeometric kernel can be written as
$$
K(x,y)=\sum_{a\in \Z'_+}\psi_a(x)\psi_a(y),
$$
which means that $K(x,y)$ is the kernel (or simply the matrix) of the spectral
projection operator in $\ell^2(\Z')$ corresponding to the positive part of the
spectrum of $D$. This formula together with a three--term recurrence relation
satisfied by the eigenfunctions $\psi_a$ implies another expression for the
kernel:
$$
K(x,y)=\frac{\sqrt{zz'\xi}}{1-\xi}\; \frac{\psi_{-\frac12}(x)\psi_{\frac12}(y)
-\psi_{\frac12}(x)\psi_{-\frac12}(y)}{x-y}\,.
$$
Thus, $K(x,y)$ is a discrete {\it integrable\/} kernel (see \cite{B} for the
definition).

\smallskip
(4) The above sum expresses the kernel as a series of products of
hypergeometric functions. On the other hand we can represent the kernel by a
double contour integral involving elementary functions only.
\smallskip

{}From these claims one can readily derive all known results concerning the
discrete hypergeometric kernel.

\subhead{The method}\endsubhead Our approach relies on the observation made in
\cite{BO2} which relates the z--measures to the Meixner orthogonal polynomials.
Assume one of the parameters $(z,z')$ is a positive integer $N=1,2,\dots$ while
the other parameter is a real number greater than $N-1$. This is a rather
special degenerate case: the weight assigned to a partition $\la$ vanishes
unless $\la_{N+1}=\la_{N+2}=\dots=0$, so that the relevant partitions $\la$
depend only on the first $N$ coordinates $\la_1,\dots,\la_N$. It turns out that
in this case the random $N$--point configuration
$\{\la_1+N-1,\,\la_2+N-2,\,\dots,\la_N\}$ on the set $\Z_+$ of nonnegative
integers is a well--known object: it is an example of a (discrete) {\it
orthogonal polynomial ensemble\/}. The orthogonal polynomial ensembles were
extensively studied in connection with random matrix theory as well as various
discrete probabilistic models, see \cite{De2}, \cite{Jo2}, \cite{Jo3},
\cite{K\"o}. In particular, it is well known that they are determinantal
processes and their correlation kernels are closely related to the
Christoffel--Darboux kernels for the corresponding family of orthogonal
polynomials (in our situation these are the classical Meixner polynomials).

The idea of our approach to the z--measures is to regard them as the result of
an {\it analytic continuation\/} of the Meixner orthogonal polynomial ensembles
in parameter $N$. In particular, our difference operator $D$ on the lattice
$\Z'$ comes from the Meixner difference operator on $\Z_+$. It is worth noting,
however, that the procedure of analytic continuation is rather delicate,
because we extrapolate from the discrete values $z=N=1,2,\dots$ to continuous
values $z\in\C$. It is this analytic continuation procedure that we regard as
the main achievement of the present paper. Even though we use it to rederive  a
known result, in a more complicated dynamical situation (see the next
paragraph) this method is crucial for obtaining new results.

Note that instead of the Meixner polynomials one could equally well use the
Krawtchouk orthogonal polynomials (see section 4).

Note also that analytic continuation of a correlation kernel off the integral
values of a parameter was used in \cite{Ni} in a very different situation. We
are grateful to the referee for this remark.

\subhead Dynamics\endsubhead The present paper can be viewed as an introduction
to our paper \cite{BO7} where the same approach is applied to studying a
dynamical model related to the z--measures. There we derive a dynamical (i.e.,
time--dependent) version of the discrete hypergeometric kernel, $K(s,x;t,y)$,
where $x$ and $y$ are, as before, points of the lattice $\Z'$ while  $s$ and
$t$ are time variables. We also evaluate the asymptotics of the kernel
$K(s,x;t,y)$ in two limit regimes. We refer to \cite{BO7} for more details.

The difference operator $D$ introduced in the present paper plays an important
role in the dynamical picture, too. We regard this operator as the key to
understanding the nature of the point processes connected to the z--measures.

\subhead Plancherel measure \endsubhead In the limit as $z$ and $z'$ go to
$\infty$ and $\xi$ goes to 0 in such a way that the product $zz'\xi$ converges
to a positive number $\theta$, the z--measure $M_{z,z',\xi}$ turns into the
poissonized Plancherel measure $M_\th$ with Poisson parameter $\theta$. Our
results about the correlation kernel have counterparts for $M_\th$, see
\cite{BOO}. A dynamical model related to $M_\th$ is studied in \cite{BO8}.

\subhead Organization of the paper\endsubhead In section 1 we recall the
definition of the z--measures and explicitly describe their relationship to the
Meixner orthogonal polynomials ensembles. In section 2 we introduce the
difference operator $D$ and we study in detail its eigenfunctions $\psi_a(x)$.
In section 3 we compute the correlation kernel. In section 4 we briefly discuss
the relationship between the z--measures and the Krawtchouk orthogonal
polynomial ensembles.

\subhead Acknowledgements \endsubhead Both authors were partially supported by
the CRDF grant RUM1-2622-ST-04. The first author (A.~B.) was also partially
supported by the NSF grant DMS-0402047.

\head 1. Z-measures
\endhead

As in Macdonald \cite{Ma} we identify partitions and Young
diagrams. By $\Y_n$ we denote the set of partitions of a natural
number $n$, or equivalently, the set of Young diagrams with $n$
boxes. By $\Y$ we denote the set of all Young diagrams, that is,
the disjoint union of the finite sets $\Y_n$, where
$n=0,1,2,\dots$ (by convention, $\Y_0$ consists of a single
element, the empty diagram $\varnothing$). Given $\la\in\Y$, let
$|\la|$ denote the number of boxes of $\la$ (so that
$\la\in\Y_{|\la|}$),  let $\ell(\la)$ be the number of nonzero
rows in $\la$ (the length of the partition), and let $\la'$
denote the transposed diagram.

By $\dim\la$ we denote the number of standard tableaux of shape $\la$. A
convenient explicit formula for $\dim\la$ is
$$
\dim\la=\frac{|\la|!}{\prod_{i=1}^N(\la_i+N-i)!}\,\prod_{1\le i<j\le
N}(\lambda_i-i-\la_j+j), \tag1.1
$$
where $N$ is an arbitrary integer $\ge \ell(\la)$ (the above
expression is stable in $N$).

We shall need the {\it generalized Pochhammer symbol\/} $(z)_\la$:
$$
(z)_\la=\prod_{i=1}^{\ell(\la)}(z-i+1)_{\la_i}\,, \qquad z\in\C,
\quad \la\in\Y,
$$
where
$$
(x)_k=x(x+1)\dots(x+k-1)=\frac{\Ga(x+k)}{\Ga(x)}
$$
is the conventional Pochhammer symbol. Note that
$$
(z)_\la=\prod_{(i,j)\in\la}(z+j-i)
$$
(product over the boxes of $\la$), which implies at once the
symmetry relation
$$
(z)_\la=(-1)^{|\la|}(-z)_{\la'}. \tag1.2
$$
Obviously, if $z=N$, where $N=1,2,\dots$, then $(z)_\la$ vanishes for all $\la$
with $\ell(\la)>N$. Likewise, if $z=-N$ then $(z)_\la$ vanishes when
$\ell(\la')=\la_1>N$.

\example{Definition 1.1} The {\it z--measure\/} with parameters $z$, $z'$, and
$\xi$ is the (complex) measure $M_{z,z',\xi}$ on the set $\Y$ which assigns to
a diagram $\la\in\Y$ the weight
$$
M_{z,z',\xi}(\la)=(1-\xi)^{zz'}\,\xi^{|\la|}\,(z)_\la(z')_\la\,
\left(\frac{\dim\la}{|\la|!}\right)^2\,. \tag1.3
$$
\endexample

The above expression makes sense for any complex $z,z'$ and any
$\xi\in\C\setminus[1,+\infty)$: Indeed, we may assume $-\pi<\arg(1-\xi)<\pi$
and then we set
$$
(1-\xi)^{zz'}=|1-\xi|^{zz'}\,e^{izz'\arg(1-\xi)}.
$$

Note that the weight is invariant under transposition $z\leftrightarrow z'$.
Note also the symmetry relation
$$
M_{-z,-z',\xi}(\la)=M_{z,z',\xi}(\la'),
$$
which readily follows from \tht{1.2}. Finally, note that the z--measures are a
particular case of the Schur measures introduced in \cite{Ok1}.

\proclaim{Proposition 1.2} If\/ $0<\xi<1$ and parameters $z,z'$ satisfy one of
the three conditions listed below, then the z--measure is a probability measure
on $\Y$. \endproclaim

The conditions are as follows.

$\bullet$ {\it Principal series\/}: The numbers $z,z'$ are not real and are
conjugate to each other.

$\bullet$ {\it Complementary series\/}: Both $z,z'$ are real and are contained
in the same open interval of the form $(m,m+1)$, where $m\in\Z$.

$\bullet$ {\it Degenerate series\/}: One of the numbers $z,z'$ (say, $z$) is a
nonzero integer while $z'$ has the same sign and, moreover, $|z'|>|z|-1$.

\demo{Proof} As follows from  \cite{BO5, \S1},  the series $\sum_\la
M_{z,z'\xi}(\la)$ absolutely converges and its sum equals 1 for any complex
$z,z'$ and any complex $\xi$ with $|\xi|<1$. Thus, it suffices to check that
the weights are nonnegative under the assumptions listed above. Since
$\xi\in(0,1)$, this means that the product $(z)_\la(z')_\la$ is nonnegative.

For the principal series, $(z)_\la$ and $(z')_\la$ are conjugate to each other
and do not vanish, and for the complementary series these are both real numbers
of the same sign. Thus, their product is always strictly positive.

Examine now the case of the degenerate series. Assume $z=N=1,2,\dots$ and
$z'>N-1$. If $\ell(\la)\le N$ then both $(z)_\la$ and $(z')_\la$ are strictly
positive, and if $\ell(\la)>N$ then $(z)_\la=0$ so that the weight vanishes.
Likewise, if $z=-N$ and $z'<-(N-1)$ then the weight is strictly positive if
$\ell(\la')=\la_1$ does exceed $N$, and vanishes otherwise. \qed
\enddemo

From now on we assume that the z--measure belongs to one of these three series
and is, therefore, a probability measure. Consequently, we may speak about {\it
random\/} Young diagrams, with reference to the z--measure.

As it is seen from the above proof, for the principal series or the
complementary series, the support of the z--measure is the whole set $\Y$,
while for the degenerate series, the support is a proper infinite subset of
$\Y$.

In the remaining part of the section we will describe the relationship between
the degenerate series and the Meixner polynomials. We start with a general
definition.

\example{Definition 1.2} Let $\X$ be a discrete subset of $\R$, finite or
countable, and let $W(x)$ be a positive function on $\X$. The $N$--point {\it
orthogonal polynomial ensemble\/} with weight function $W$ is the random
$N$--point configuration in $\X$ such that the probability of a particular
configuration $x_1>\dots>x_N$, where $x_1,\dots,x_N\in\X$, is given by
$$
\Prob\{x_1,\dots,x_N\}=\const_N\, \prod_{i=1}^N W(x_i)\, \prod_{1\le i<j\le
N}(x_i-x_j)^2. \tag1.4
$$
\endexample

Here we assume that the cardinality of $\X$ is no less than $N$ and that
$$
\const_N^{-1}:=\sum_{x_1>\dots>x_N} \left\{\prod_{i=1}^N W(x_i)\, \prod_{1\le
i<j\le N}(x_i-x_j)^2\right\}<+\infty.
$$
For finite $\X$, this condition is trivial, and for infinite $\X$, it just
means that the weight function $W$ has at least $N-1$ finite first moments.

The term ``orthogonal polynomial ensemble'' is related to the following
well--known fact. Let $P_0=1, P_1,P_2,\dots$ be the orthogonal polynomials with
weight function $W$, the result of Gram--Schmidt orthogonalization of
$1,x,x^2,\dots$ in the weighted $\ell^2$ space $\ell^2(\X,W)$. Denote by
$K_N(x,y)$ the $N$th Christoffel--Darboux kernel multiplied by
$\sqrt{W(x)W(y)}$:
$$
K_N(x,y)=\sqrt{W(x)W(y)}\sum_{i=0}^{N-1}\frac{P_i(x)P_i(y)}{\Vert P_i\Vert^2},
$$
where the norm refers to the weighted $\ell^2$ Hilbert space $\ell^2(\X,W)$.
Note that the kernel $K_N$ corresponds to the projection operator in
$\ell^2(\X,W)$ whose range is the linear span of $1,x,\dots,x^{N-1}$. Then we
have

\proclaim{Proposition 1.3} The probability that the random $N$--point
configuration, as specified in Definition 1.2, contains a given $n$-point set
$\{y_1,\dots,y_n\}\subset\X$ equals the determinant of the $n\times n$ matrix
$[K_N(y_i,y_j)]$.
\endproclaim

For a proof, see, e.g., \cite{De2}, \cite{K\"o, Lemma 2.8}. Note that the
determinant automatically vanishes if $n>N$, because the kernel has rank $N$.

We will be dealing with a concrete example of the weight function. This is the
{\it Meixner weight function\/}, which is defined on the set
$\X=\Z_+:=\{0,1,\dots\}$, depends on parameters $\be>0$ and $\xi\in(0,1)$, and
is given by
$$
W^\m_{\be,\xi}(\ti x)= \frac{(\be)_{\ti x}\xi^{\ti x}}{\ti x!}
=\frac{\Ga(\be+\ti x)\xi^{\ti x}}{\Ga(\be) \ti x!}\,, \qquad  \ti x\in\Z_+
$$
(we denote a point of $\Z_+$ by $\ti x$ instead of $x$ because this notation is
used below in \S\S2--3).

For $N=1,2,\dots$, let $\Y(N)\subset\Y$ denote the set of diagrams $\la$ with
$\ell(\la)\le N$. The following correspondence is a bijection between diagrams
$\la\in\Y(N)$ and $N$--point configurations on $\Z_+$:
$$
\la\,\longleftrightarrow\, \{\ti x_1,\dots,\ti x_N\}, \qquad \ti x_i=\la_i+N-i
\quad (i=1,\dots,N). \tag1.5
$$

The next fact was pointed out in \cite{BO2}:

\proclaim{Proposition 1.4} Under correspondence\/ \tht{1.5}, the z--measure of
the degenerate series with parameters\/ $(z=N,\, z'=N+\be-1, \,\xi)$, where
$\be>0$ and $\xi\in(0,1)$,  turns into the $N$--point Meixner orthogonal
polynomial ensemble with parameters $(\be,\xi)$.
\endproclaim

\demo{Proof} It suffices to check that if $\ell(\la)\le N$ and $\{\ti
x_1,\dots,\ti x_N\}$ is given by \tht{1.5} then the right--hand side of
\tht{1.3} can be written as the right--hand side of \tht{1.4} with the Meixner
weight function.

By virtue of \tht{1.1},
$$
\left(\frac{\dim\la}{|\la|!}\right)^2=\frac{\prod_{i<j}(\ti x_i-\ti
x_j)^2}{\prod_i (\ti x_i!)^2}\,.
$$
Next, with $z=N$ and $z'=N+\be-1$ we have
$$
(z)_\la=\frac{\prod_{i=1}^N \ti x_i!}{\prod_{i=1}^N (N-i)!}\,, \qquad
(z')_\la=\frac{\prod_{i=1}^N \Ga(\be+\ti x_i)}{\prod_{i=1}^N \Ga(\be+N-i)}\,,
$$
and
$$
\xi^{|\la|}=\xi^{-N(N-1)/2}\,\prod_{i=1}^N\xi^{\ti x_i}.
$$
Combining these formulas we get \tht{1.4} with $W=W^\m_{\be,\xi}$ and
$$
\const_N=(1-\xi)^{N(N+\be-1)}\xi^{-N(N-1)/2}\prod_{i=1}^N
\frac{\Ga(\be)}{\Ga(N-i+1)\Ga(N-i+\be)}\,.
$$
\qed
\enddemo

\example{Remark 1.5} Let $\la$ be the random Young diagram distributed
according to a z--measure $M_{z,z',\xi}$. Then the number of boxes $|\la|$ has
the negative binomial distribution on $\Z_+$ with parameters $zz'$ and $\xi$:
$$
\multline \Prob\{|\la|=n\}=\pi_{zz',\xi}(n)\\
:=(1-\xi)^{zz'}\cdot\xi^n\cdot \frac{zz'(zz'+1)\dots(zz'+n-1)}{n!}\,.
n=0,1,2,\dots,
\endmultline
$$
Conditioned on $|\la|=n$, the distribution of $\la$ is a probability measure
$M^{(n)}_{z,z'}$ on $\Y_n$ which does not depend on $\xi$:
$$
M^{(n)}_{z,z'}(\la)=\frac{(z)_\la(z')_\la}{zz'(zz'+1)\dots(zz'+n-1)}
\cdot\frac{(\dim\la)^2}{n!}\,, \qquad \la\in\Y_n
$$
(recall that $\Y_n$ is the set of diagrams with $n$ boxes). This means that the
z--measure $M_{z,z',\xi}$ is the mixture of the probability measures
$M^{(n)}_{z,z'}$ with varying index $n\in\Z_+$ by means of the negative
binomial distribution $\pi_{zz',\xi}$, see \cite{BO2}, \cite{BO4}. For
applications to harmonic analysis on the infinite symmetric group one needs the
measures $M^{(n)}_{z,z'}$ and their scaling limits as $n\to\infty$, but it
turns out that the ``mixed'' measures $M_{z,z',\xi}$ have much better
properties, and the large $n$ limit can be replaced, to a certain extent, by
the $\xi\to1$ limit transition. This was the starting point of our paper
\cite{BO2}. In the present paper we are dealing with the ``mixed'' measures
only.

\endexample

\head 2. A basis in the $\ell^2$ space on the lattice and the
Meixner polynomials \endhead

In this section we examine a nice orthonormal basis in the
$\ell^2$ space on the 1--dimensional lattice. The elements of
this basis are eigenfunctions of a second order difference
operator. They can be obtained from the classical Meixner
polynomials via analytic continuation with respect to
parameters.

Throughout the section we will assume (unless otherwise stated) that parameters
$z,z'$ are in the principal series or in the complementary series but not in
the degenerate series. In particular, $z,z'$ are not integers.

Consider the lattice of (proper) half--integers
$$
\Z'=\Z+\tfrac12=\{\dots,-\tfrac52,-\tfrac32,-\tfrac12,
\,\tfrac12,\,\tfrac32,\,\tfrac52,\dots\}.
$$
Elements of $\Z'$ will be denoted by letters $x,y$.

We introduce a family of functions on $\Z'$ depending on a parameter $a\in\Z'$
and also on our basic parameters $z,z',\xi$:
$$
\gather
\psi_a(x;z,z',\xi)=\left(\frac{\Ga(x+z+\tfrac12)\Ga(x+z'+\tfrac12)}
{\Ga(z-a+\tfrac12)\Ga(z'-a+\tfrac12)}\right)^{\frac12}
\xi^{\frac12(x+a)}(1-\xi)^{\tfrac12(z+z')-a}\\
\times\frac{F(-z+a+\tfrac12,-z'+a+\tfrac12;x+a+1;\w)}{\Ga(x+a+1)}\,, \qquad
x\in\Z', \tag2.1
\endgather
$$
where $F(A,B;C;w)$ is the Gauss hypergeometric function.

Let us explain why this expression makes sense. Since, by convention,
parameters $z,z'$ do not take integral values, $\Ga(x+z+\tfrac12)$ and
$\Ga(x+z'+\tfrac12)$ have no singularities for $x\in\Z'$. Moreover, the
assumptions on $(z,z')$ imply that
$$
\Ga(x+z+\tfrac12)\Ga(x+z'+\tfrac12)>0, \qquad
\Ga(z-a+\tfrac12)\Ga(z'-a+\tfrac12)>0,
$$
so that we can take the positive value of the square root in
\tht{2.1}. Next, since $\xi\in(0,1)$, we have $\xi/(\xi-1)<0$,
and as is well known, the function $w\to F(A,B;C,w)$ is well
defined on the negative semi--axis $w<0$. Finally, although
$F(A,B;C,w)$ is not defined at $C=0,-1,-2,\dots$, the  ratio
$F(A,B;C,w)/\Ga(C)$ is well defined for all $C\in\C$.

Note also that the functions $\psi_a(x;z,z',\xi)$ are
real--valued. Their origin will be explained below.

Further, we introduce a second order difference operator
$D(z,z',\xi)$ on the lattice $\Z'$, depending on parameters
$z,z',\xi$ and acting on functions $f(x)$ (where $x$ ranges over
$\Z'$) as follows
$$
\gather D(z,z',\xi)f(x)=\sqrt{\xi(z+x+\tfrac12)(z'+x+\tfrac12)}\,f(x+1)\\
+\sqrt{\xi(z+x-\tfrac12)(z'+x-\tfrac12)}\,f(x-1)
-(x+\xi(z+z'+x))\,f(x).
\endgather
$$
Note that $D(z,z',\xi)$ is a symmetric operator in
$\ell^2(\Z')$.

\proclaim{Proposition 2.1} The functions $\psi_a(x;z,z',\xi)$,
where $a$ ranges over $\Z'$, are eigenfunctions of the operator
$D(z,z',\xi)$,
$$
D(z,z',\xi)\psi_a(x;z,z',\xi)=a(1-\xi)\psi_a(x;z,z',\xi).
\tag2.2
$$
\endproclaim

\demo{Proof} This equation can be verified using the relation
$$
\gather
w(C-A)(C-B)F(A,B;C+1;w)-(1-w)C(C-1)F(A,B;C-1;w)\\
+C[C-1-(2C-A-B-1)w]F(A,B;C;w)=0
\endgather
$$
for the Gauss hypergeometric function, see, e.g., \cite{Er, 2.8
(45)}. \qed
\enddemo

The next lemma  provides us a convenient integral representation for the
functions $\psi_a$.

\proclaim{Lemma 2.2} For any $A,B\in\C$, $M\in\Z$, and
$\xi\in(0,1)$ we have
$$
\gathered
\frac{F(A,B;M+1;\w)}{\Ga(M+1)}=\frac{\Ga(-A+1)\xi^{-M/2}(1-\xi)^B}
{\Ga(-A+M+1)}\\
\times\, \frac1{2\pi i}\int\limits_{\{\om\}}
(1-\sqrt\xi\om)^{A-1}\left(1-\frac{\sqrt\xi}{\om}\right)^{-B}
\om^{-M}\,\frac{d\om}{\om}\,.
\endgathered \tag2.3
$$
Here $\xi\in(0,1)$ and $\{\om\}$ is an arbitrary simple contour
which goes around the points 0 and $\sqrt\xi$ in the positive
direction leaving $1/\sqrt\xi$ outside.
\endproclaim

\demo{Comments} 1. The branch of the function
$(1-\sqrt\xi\om)^{A-1}$ is specified by the convention that the
argument of $1-\sqrt\xi\om$ equals 0 for real negative values of
$\om$, and the same convention is used for the function
$\left(1-\frac{\sqrt\xi}{\om}\right)^{-B}$.

2. Like the Euler integral formula, formula \tht{2.3} does not
make evident the symmetry $A\leftrightarrow B$.

3. The right--hand side of formula \tht{2.3} makes sense for
$A=1,2,\dots$, when $\Ga(-A+1)$ has a singularity. Then the
whole expression can be understood, e.g., as the limit value as
$A$ approaches one of the points 1,2,\dots .
\enddemo

\demo{Proof} Since both sides of \tht{2.3} are real--analytic
functions of $\xi$ we may assume that $\xi$ is small enough.
Then we may apply the binomial formula which gives
$$
\gather
\xi^{-M/2}(1-\sqrt\xi\om)^{A-1}\left(1-\frac{\sqrt\xi}{\om}\right)^{-B}
\om^{-M}\\
= \sum_{k=0}^\infty\sum_{l=0}^\infty \frac{(-A+1)_k
(B)_l}{k!\,l!}\,\xi^{(k+l-M)/2}\,\om^{k-l-M}\,.
\endgather
$$
After integration only the terms with $k=l+M$ survive. It
follows that the right--hand side of \tht{2.3} is equal to
$$
\frac{(1-\xi)^B}{\Ga(-A+M+1)}\,\sum_{l\ge\max(0,-M)}
\frac{\Ga(-A+M+1+l)\,(B)_l} {\Ga(l+M+1)l!}\, \xi^l.
$$
We may replace the inequality $l\ge\max(0,-M)$ simply by $l\ge0$
because for negative integral values of $M$ (when we have to
start summation from $l=-M$), the terms with $l=0,\dots,-M-1$
automatically vanish due to the factor $\Ga(l+M+1)$ in the
denominator. Consequently, our expression is equal to
$$
\frac{(1-\xi)^B
F(-A+1+M,B;M+1;\xi)}{\Ga(M+1)}=\frac{F(A,B;M+1;\w)}{\Ga(M+1)}\,,
$$
where we used the transformation formula \cite{Er, 2.9 (4)}. \qed
\enddemo

\proclaim{Proposition 2.3} We have the following integral
representations
$$
\multline \psi_{a}(x;z,z',\xi)\\
=\left(\frac{\Ga(x+z+\tfrac12)\Ga(x+z'+\frac12)}
{\Ga(z-a+\tfrac12)\Ga(z'-a+\frac12)} \right)^\frac
12\, \frac{\Ga(z'-a+\tfrac12)}{\Ga(z'+ x+\tfrac12)} \,(1-\xi)^{\frac{z'-z+1}2}\\
\times\frac1{2\pi
i}\,\oint\limits_{\{\om\}}\left(1-\sqrt{\xi}\om\right)^{-z'+a-\tfrac12}
\left(1-\frac{\sqrt{\xi}}{\om}\right)^{z-a-\tfrac12}\om^{-x-a}
\,\frac{d\om}{\om}
\endmultline\tag2.4
$$
and
$$
\multline \psi_a(x;z,z',\xi)\psi_a(y;z,z',\xi)=\varphi_{z,z'}(x,y)\\
\times\,\frac{1-\xi}{(2\pi i)^2}
\oint\limits_{\{\om_1\}}\oint\limits_{\{\om_2\}}
\left(1-\sqrt{\xi}\om_1\right)^{-z'+a-\tfrac12}
\left(1-\frac{\sqrt{\xi}}{\om_1}\right)^{z-a-\tfrac12}\\
\times\left(1-\sqrt{\xi}\om_2\right)^{-z+a-\tfrac12}
\left(1-\frac{\sqrt{\xi}}{\om_2}\right)^{z'-a-\tfrac12}
\om_1^{-x-a}\om_2^{-y-a}
\,\frac{d\om_1}{\om_1}\frac{d\om_2}{\om_2}
\endmultline \tag2.5
$$
where
$$
\varphi_{z,z'}(x,y)=\frac{\sqrt{\Ga(x+z+\tfrac12)\Ga(x+z'+\tfrac12)
\Ga(y+z+\tfrac12)\Ga(y+z'+\tfrac12)}}
{\Ga(x+z'+\tfrac12)\Ga(y+z+\tfrac12)} \tag2.6
$$
Here each contour is an arbitrary simple loop, oriented in
positive direction, surrounding the points 0 and $\sqrt\xi$, and
leaving $1/\sqrt\xi$ outside. We also use the convention about
the choice of argument as in Comment 1 to Lemma 2.2.
\endproclaim

\demo{Proof} Indeed, \tht{2.4} immediately follows from
\tht{2.1} and \tht{2.3}. To prove \tht{2.5} we multiply out the
integral representation \tht{2.4} for the first function and the
same representation for the second function, but with $z$ and
$z'$ interchanged. The transposition $z\leftrightarrow z'$ in
\tht{2.4} is justified by the fact the initial formula \tht{2.1}
is symmetric with respect to $z\leftrightarrow z'$. As a result
of this trick the gamma prefactors involving $a$ are completely
cancelled out, and we obtain \tht{2.5}\qed
\enddemo

\proclaim{Proposition 2.4} The functions
$\psi_a=\psi_a(x;z,z',\xi)$, where $a$ ranges over $\Z'$, form
an orthonormal basis in the Hilbert space $\ell^2(\Z')$.
\endproclaim

\demo{Proof} {}From \tht{2.4} it is not difficult to see that the function
$\psi_a(x;z,z',\xi)$ has exponential decay as $x\to\pm\infty$. Indeed,
depending on whether $x$ goes to $+\infty$ or $-\infty$ we arrange the contour
in such a way that $|\om|>1$ or $|\om|<1$, respectively.

In particular, $\psi_a(x;z,z',\xi)$ is square integrable. Since $\psi_a$ is an
eigenfunction of a symmetric difference operator whose coefficients have linear
growth at $\pm\infty$, and since to different indices $a$ correspond different
eigenvalues, we conclude that these functions are pairwise orthogonal in
$\ell^2(\Z')$.

Let us show that $\Vert\psi_a\Vert^2=1$. Take \tht{2.5} with $x=y$. Then the
whole expression simplifies because \tht{2.6} turns into 1. Next, in the double
contour integral, we replace the variable $\om_2$ by its inverse. We obtain
$$
\gathered (\psi_{a}(x;z,z',\xi))^2=\frac{1-\xi}{(2\pi i)^2}
\oint\oint\left(1-\sqrt{\xi}\om_1\right)^{-z'+a-\tfrac12}
\left(1-\sqrt{\xi}\,\om_1^{-1}\right)^{z-a-\tfrac12}\\
\times\left(1-\sqrt{\xi}{\om_2}^{-1}\right)^{-z+a-\tfrac12}
\left(1-\sqrt{\xi}\,\om_2\right)^{z'-a-\tfrac12}
\left(\frac{\om_1}{\om_2}\right)^{-x-a}
\,\frac{d\om_1}{\om_1}\frac{d\om_2}{\om_2}
\endgathered
$$
To evaluate the squared norm we have to sum this expression over
$x\in\Z'$. We split the sum into two parts according to the
splitting $\Z'=\Z'_-\cup\Z'_+$. We take as the contours
concentric circles such that $|\om_1|<|\om_2|$ in the sum over
$\Z'_-$, and $|\om_1|>|\om_2|$ in the sum over $\Z'_+$. This
gives us
$$
\sum\limits_{x\in\Z'}(\psi_{a}(x;z,z',\xi))^2=
\underset{|\om_1|<|\om_2|}\to{\oint\oint}\frac{F(\om_1,\om_2)}{\om_2-\om_1}\,
\frac{d\om_1}{\om_1}\frac{d\om_2}{\om_2} +
\underset{|\om_1|>|\om_2|}\to{\oint\oint}\frac{F(\om_1,\om_2)}{\om_1-\om_2}\,
\frac{d\om_1}{\om_1}\frac{d\om_2}{\om_2}
$$
with
$$
\gather F(\om_1,\om_2)=\frac{1-\xi}{(2\pi
i)^2}\,\left(1-\sqrt{\xi}\om_1\right)^{-z'+a-\frac12}
\left(1-\sqrt{\xi}\,\om_1^{-1}\right)^{z-a-\frac12}\\
\times\left(1-\sqrt{\xi}{\om_2}^{-1}\right)^{-z+a-\frac12}
\left(1-\sqrt{\xi}\,\om_2\right)^{z'-a-\frac12}\,
\om_1^{\frac12-a}\om_2^{\frac12+a}.
\endgather
$$
Recall that both contours go in positive direction.

Let us transform the second double--contour integral: keeping the second
contour fixed we move the first contour inside the second contour. Then we
obtain a double--contour integral which cancels the first double--contour
integral, plus a single--contour integral arising from the residue of the
function $\om_1\to(\om_1-\om_2)^{-1}$:
$$
\frac{1-\xi}{2\pi i}\oint
F(\om,\om)\frac{d\om}{\om^2}=\frac{1-\xi}{2\pi
i}\oint\frac{d\om}{(1-\sqrt\xi\om)(\om-\sqrt\xi)}=1.
$$

Thus, we have shown that the functions $\psi_a$ form an
orthonormal family in $\ell^2(\Z')$, and it remains to prove
that this family is complete. For $x\in\Z'$, let $\de_x$ stand
for the delta function at $x$. Since the functions $\de_x$ form
an orthonormal basis, it suffices to check that
$$
\sum_{a\in\Z'}\left((\de_x,\psi_a)_{\ell^2(\Z')}\right)^2 =
\sum_{a\in\Z'}(\psi_a(x;z,z',\xi))^2=1, \qquad \forall x\in\Z'.
$$
But this follows from the previous claim and  the symmetry
$a\leftrightarrow x$ established in the next proposition. \qed
\enddemo

\proclaim{Proposition 2.5} The following symmetry relation holds
$$
\psi_a(x;z,z',\xi)=\psi_x(a;-z,-z',\xi).
$$
\endproclaim

\demo{Proof} Using the classical formula
$$
\Ga(A+\tfrac12)\Ga(A-\tfrac12)=\frac{\pi}{\cos(\pi A)}
$$
and the fact that both $x+\tfrac12$ and $a+\tfrac12$ are
integers we check that
$$
\frac{\Ga(z+x+\tfrac12)\Ga(z'+x+\tfrac12)}
{\Ga(z-a+\tfrac12)\Ga(z'-a+\tfrac12)}
=\frac{\Ga(-z+a+\tfrac12)\Ga(-z'+a+\tfrac12)}
{\Ga(-z-x+\tfrac12)\Ga(-z'-x+\tfrac12)}\,.
$$
Applying this to \tht{2.1} and using another classical formula,
$$
F(A,B;C;w)=(1-w)^{C-A-B}F(C-A,C-B;C;w),
$$
see \cite{Er, (2.9.2)}, we get the required relation.

Another way to prove the proposition is to make a change of the variable in
integral \tht{2.4}:
$$
\om\mapsto \om'=\frac{\om-\sqrt\xi}{\sqrt\xi\om-1}\,.
$$
This is an involutive transformation such that
$0\leftrightarrow\sqrt\xi$ and $\infty\leftrightarrow
1/\sqrt\xi$. As is readily verified, it leads to transformation
$(a,x,z,z')\to(x,a,-z,-z')$. \qed
\enddemo

\proclaim{Corollary 2.6} The functions
$\psi_a=\psi_a(x;z,z',\xi)$ satisfy the following three--term
relation
$$
\gather (1-\xi)x \psi_a
=\sqrt{\xi(z-a+\tfrac12)(z'-a+\tfrac12)}\,\psi_{a-1}\\
+\sqrt{\xi(z-a-\tfrac12)(z'-a-\tfrac12)}\,\psi_{a+1}
+(-a+\xi(z+z'-a))\,\psi_a.
\endgather
$$
\endproclaim

\demo{Proof} Under symmetry $x\leftrightarrow a$ (Proposition
2.4), this turns into the formula stated in Proposition 2.1. Of
course, a direct verification is also possible. \qed
\enddemo

The formulas of Proposition 2.1 and Corollary 2.6 show that the
functions $\psi_a(x;z,z',\xi)$ possess the {\it bispectrality\/}
property in the sense of \cite{Gr}.

\proclaim{Proposition 2.7} One more symmetry relation holds:
$$
\psi_a(x;z,z',\xi)=(-1)^{x+a}\psi_{-a}(-x;-z,-z',\xi), \qquad
x,a\in\Z'.
$$
\endproclaim

\demo{Proof} This follows from the relation
$$
\gather\frac{F(A,B;C;w)}{\Ga(C)} =w^{1-C}\,
\frac{\Ga(A-C+1)\Ga(B-C+1)}{\Ga(A)\Ga(B)}\\
\times\frac{F(A-C+1,B-C+1;2-C;w)}{\Ga(2-C)}\,, \qquad C\in\Z,
\endgather
$$
see \cite{Er, 2.8 (19)}. Another way is to make a change of the
variable, $\om\mapsto1/\om$, in integral \tht{2.4}. \qed
\enddemo

In the remaining part of the section we will explain how the
functions $\psi_a$ are related to the Meixner polynomials.

Let $\Z_+=\{0,1,2,\dots\}$. To denote points of $\Z_+$ we will use now the
symbols $\ti x, \ti y$, because the letters $x,y$ were already employed to
denote points of $\Z'$. Recall that the Meixner polynomials are the orthogonal
polynomials with respect to the weight function
$$
W^\m_{\be,\xi}(\ti x)= \frac{(\be)_{\ti x}\xi^{\ti x}}{\ti x!}
=\frac{\Ga(\be+\ti x)\xi^{\ti x}}{\Ga(\be)\ti x!}\,,  \qquad \ti x\in\Z_+\,,
\tag2.7
$$
on $\Z_+$, where $\be>0$ and, as before, $\xi\in(0,1)$. Our notation for these
polynomials is $\Me_n(\ti x;\be,\xi)$. We use the same normalization of the
polynomials as in the handbook \cite{KS} (note that in \cite{KS}, our parameter
$\xi$ is denoted as $c$).

Set
$$
\wt\Me_n(\ti x;\be,\xi)=(-1)^n\,\frac{\Me_n(\ti
x;\be,\xi)}{\Vert\Me_n(\,\cdot\,;\be,\xi)\Vert}\,\sqrt{W^\m_{\be,\xi}(\ti x)},
\qquad \ti x\in\Z_+\,,\tag2.8
$$
where
$$
\Vert\Me_n(\,\cdot\,;\be,\xi)\Vert^2 =\sum_{\ti x=0}^\infty \Me_n^2(\ti
x;\be,\xi)W^\m_{\be,\xi}(\ti x).
$$
The factor $(-1)^n$ is introduced for convenience: it will
compensate the same factor in formula \tht{2.10} below.

\proclaim{Proposition 2.8} Drop the assumption that $(z,z')$ is not in the
degenerate series, and assume, just on the contrary, that $z=N$ and
$z'=N+\be-1$, where $N=1,2,\dots$ and $\be>0$. Then expression \tht{2.1} for
the functions $\psi_a(x;z,z',\xi)$ still makes sense provided that the numbers
$$
\ti x:=x+N-\tfrac12\,, \qquad n:=N-a-\tfrac12 \tag2.9
$$
are in $\Z_+$, and in this notation we  have
$$
\psi_a(x;z,z',\xi)=\wt\Me_n(\ti x;\be,\xi).
$$
\endproclaim

\demo{Proof} We start with the expression of the Meixner
polynomials through the Gauss hypergeometric function (see
\cite{KS, \tht{1.9.1}}):
$$
\Me_n(\ti x;\be,\xi)=F(-n,-\ti x;\be;\tfrac{\xi-1}{\xi}).
$$
Applying the transformation
$$
F(-n,b;\be;w)=\frac{\Ga(1-b)\Ga(\be)w^n}{\Ga(\be+n)}\,
\frac{F(-n,1-\be-n;1-b-n;w^{-1})}{\Ga(1-b-n)}\,, \quad n=0,1,\dots,
$$
we obtain
$$
\multline \Me_n(\ti x;\be,\xi) \\
=\frac{(-1)^n\Ga(\ti x+1)\Ga(\be)}{\Ga(\be+n)}\, \left(
\frac{1-\xi}{\xi}\right)^n\, \frac{F(-n,-\be-n+1;\ti x+1-n;\w)}{\Ga(\ti
x+1-n)}\,.
\endmultline\tag2.10
$$
Although the first expression for the polynomials looks simpler
than the second one, it turns out that only the second
expression is suitable for our purposes. Note that (see
\cite{KS, \tht{1.9.2}})
$$
\Vert \Me_n(\,\cdot\,;\be,\xi)\Vert^{-2}
=\frac{\xi^n(1-\xi)^{\be}\Ga(\be+n)}{\Ga(\be)\Ga(n+1)}\,.
$$

{}From the last two formulas and the definition of $\wt\Me_n$ we obtain
$$
\gather \wt\Me_n(\ti x;\be,\xi)=\sqrt{\frac{\Ga(\ti x+1)\Ga(\ti
x+\be)}{\Ga(n+1)\Ga(n+\be)}}\,\xi^{(\ti
x-n)/2}(1-\xi)^{(\be+2n)/2}\\
\times \frac{F(-n,-\be-n+1;\ti x+1-n;\w)}{\Ga(\ti x+1-n)}\,.
\endgather
$$
Comparing this with \tht{2.1} and taking into account \tht{2.9}
we get the required equality. \qed
\enddemo

Thus, our functions $\psi_a$ can be obtained from the Meixner
polynomials by the following procedure:

$\bullet$ We replace the initial polynomials $\Me_n$ by the functions
$\wt\Me_n$. This step is quite clear: as a result we get functions which form
an orthonormal basis in the $\ell^2$ space on $\Z_+$ with respect to the weight
function 1.

$\bullet$ Next, we make a change of the argument. Namely, we
introduce an additional parameter $N=1,2,\dots$ and we set
$x=\ti x-N+\tfrac12$. Then we get orthogonal functions on the
subset
$$
\{-N+\tfrac12, -N+\tfrac32, -N+\tfrac52, \dots\}\subset\Z',
$$
which exhausts the whole $\Z'$ in the limit as $N$ goes to
infinity.

$\bullet$ Then we also need a change of the index. Namely,
instead of $n$ we have to take $a=N-n-\tfrac12$. We cannot give
a conceptual explanation of this transformation, it is dictated
by the formulas. Again,  the range of the possible values for
$a$ becomes larger together with $N$, and in the limit as
$N\to+\infty$ we get the whole lattice $\Z'$.

$\bullet$ Finally, we make a (formal) analytic continuation in parameters $N$
and $\be$, using an appropriate analytic expression for the Meixner polynomials
(namely, \tht{2.10}).

We hope that this detailed explanation will help the reader to perceive the
analytic continuation arguments in section 3.

Of course, instead of the lattice $\Z'$ we could equally well
deal with the lattice $\Z$, and then numerous ``$\frac12$''
would disappear. However, dealing with the lattice $\Z'$ makes
main formulas more symmetric.

\example{Remark 2.9} Note that the difference equation of Proposition 2.1 can
be obtained via the procedure described above from the classical difference
equation satisfied by the Meixner polynomials. This is precisely the way how we
have obtained the difference operator $D$. Likewise, the three--term relation
of Corollary 2.6 precisely corresponds to the classical three--term relation
for the Meixner polynomials.
\endexample

\head 3. The discrete hypergeometric kernel
\endhead

Let $\X$ be a countable set. By a {\it point configuration\/} in
$\X$ we mean any subset $X\subseteq\X$. Let $\Conf(\X)$ be the
set of all point configurations; this is a compact space. Assume
we are given a probability measure on $\Conf(\X)$ so that we can
speak about the {\it random\/} point configuration in $\X$. The
$n$th {\it correlation function\/} of our probability measure
(where $n=1,2,\dots$) is defined by
$$
\rho_n(x_1,\dots,x_n)=\Prob\{\text{the random configuration
contains $x_1,\dots,x_n$}\},
$$
where $x_1,\dots,x_n$ are pairwise distinct points in $\X$. The
collection of all correlation functions determines the initial
probability measure uniquely.

We say that our probability measure is {\it determinantal\/} if
there exists a function $K(x,y)$ on $\X\times\X$ such that
$$
\rho_n(x_1,\dots,x_n)=\det\left[K(x_i,x_j)\right]_{i,j=1}^n\,,
\qquad n=1,2,\dots \tag3.1
$$
It is worth noting that if such a function $K(x,y)$ exists, then
it is not unique. Indeed, any ``gauge transformation'' of the
form
$$
K(x,y)\to \frac{f(x)}{f(y)}\,K(x,y), \tag3.2
$$
where $f$ is a nonvanishing function on $\X$, does not affect
the determinants in the right--hand side of \tht{3.1}.

Any function $K(x,y)$ satisfying \tht{3.1} will be called a {\it correlation
kernel\/} of the initial determinantal measure. Two kernels giving the same
system of correlation functions will be called {\it equivalent}.

As in \S2, we are dealing with the lattice $\Z'$ of (proper)
half--integers. We split it into two parts,
$\Z'=\Z'_-\cup\Z'_+$, where $\Z'_-$ consists of all negative
half--integers and $\Z'_+$ consists of all positive
half--integers. For an arbitrary $\la\in\Y$ we set
$$
\unX(\la)=\{\la_i-i+\tfrac12\mid i=1,2,\dots\}\subset\Z'.
$$
For instance, $\unX(\varnothing)=\Z'_-$. The set $\unX(\la)$ is sometimes
called the {\it Maya diagram\/} of $\la$, see, e.g. \cite{MJD}.

The correspondence $\la\mapsto \unX(\la)$ is a bijection between the Young
diagrams $\la$ and those (infinite) subsets $X\subset\Z'$ for which the
symmetric difference $X\triangle\,\Z'_-$ is a finite set with equally many
points in $\Z'_+$ and $\Z'_-$. Note that
$$
\unX(\la')=-(\Z'\setminus \unX(\la)).
$$

Using the correspondence  $\la\mapsto \unX(\la)$ we can interpret any
probability measure $M$ on $\Y$  as a probability measure on $\Conf(\Z')$. This
makes it possible to speak about the correlation functions of $M$. Our goal is
to compute them explicitly for the z--measures.

Now we can state the main results of the paper.

\proclaim{Theorem 3.1} Under the above correspondence between Young diagrams
and Maya diagrams, any z--measure determines a determinantal measure on
$\Conf(\Z')$.
\endproclaim

\proclaim{Theorem 3.2} The correlation kernel of any z--measure $M_{z,z',\xi}$
from the principal or complementary series can be written in the form
$$
\unK_{z,z',\xi}(x,y)=\sum_{a\in\Z'_+}\psi_a(x;z,z',\xi)\psi_a(y;z,z',\xi),
\qquad x,y\in\Z', \tag3.3
$$
where the functions $\psi_a$ are defined in \tht{2.1}.
\endproclaim

Note that the series in the right--hand side is absolutely
convergent. Indeed, since $\{\psi_a\}$ is an orthonormal basis
in $\ell^2(\Z')$ (Proposition 2.4), this follows from the fact
that the series can be written as
$$
\sum_{a\in\Z'_+}(\de_x,\psi_a)(\psi_a,\de_y),
$$
where $\de_x$ stands for the delta--function at point $x$ on the
lattice $\Z'$, and $(\,\cdot\,,\,\cdot\,)$ denotes the inner
product in $\ell^2(\Z')$.

Formula \tht{3.3} simply means that $\unK_{z,z',\xi}(x,y)$ is
the matrix of the orthogonal projection operator in
$\ell^2(\Z')$ whose range is the subspace spanned by the basis
vectors $\psi_a$ with index $a\in\Z'_+\subset\Z'$.

\proclaim{Theorem 3.3} The correlation kernel \tht{3.3} can also
be written in the form
$$
\unK_{z,z',\xi}(x,y)=\varphi_{z,z'}(x,y)\,\hatK_{z,z',\xi}(x,y)\tag3.4
$$
where, as in \tht{2.6},
$$
\varphi_{z,z'}(x,y)=\frac{\sqrt{\Ga(x+z+\tfrac12)\Ga(x+z'+\tfrac12)
\Ga(y+z+\tfrac12)\Ga(y+z'+\tfrac12)}}
{\Ga(x+z'+\tfrac12)\Ga(y+z+\tfrac12)} \tag3.5
$$
and
$$
\multline \hatK_{z,z',\xi}(x,y)\\=\frac1{(2\pi
i)^2}\oint\limits_{\{\om_1\}}\oint\limits_{\{\om_2\}}
\dfrac{(1-\sqrt\xi\om_1)^{-z'} \left(1-\dfrac{\sqrt\xi}{\om_1}\right)^{z}
(1-\sqrt\xi\om_2)^{-z}
\left(1-\dfrac{\sqrt\xi}{\om_2}\right)^{z'}}{\om_1\om_2-1}\\
\times \,\om_1^{-x-\tfrac12}\om_2^{-y-\tfrac12}d\om_1\,d\om_2
\endmultline\tag3.6
$$
where $\{\om_1\}$ and $\{\om_2\}$ are arbitrary simple contours
satisfying the following three conditions{\/\rm:}

$\bullet$ both contours go around 0 in positive
direction{\/\rm;}

$\bullet$ the point $\xi^{1/2}$ is in the interior of each of the contours
while the point $\xi^{-1/2}$ lies outside the contours{\rm;}

$\bullet$ the contour $\{\om_1^{-1}\}$ is contained in the
interior of the contour $\{\om_2\}$ {\rm(}equivalently,
$\{\om_2^{-1}\}$ is contained in the interior of
$\{\om_1\}${\rm)}.

The kernels $\unK_{z,z',\xi}(x,y)$ and $\hatK_{z,z,\xi}(x,y)$
are equivalent. Namely, they are related by a ``gauge
transformation'',
$$
\hatK_{z,z',\xi}(x,y)
=\frac{f_{z,z'}(x)}{f_{z,z'}(y)}\,\unK_{z,z',\xi}(x,y), \qquad
x,y\in\Z',
$$
where
$$
f_{z,z'}(x)=\frac{\Ga(x+z'+\tfrac12)}
{\sqrt{\Ga(x+z+\tfrac12)\Ga(x+z'+\tfrac12)}}\tag3.7
$$
The kernel $\hatK_{z,z',\xi}(x,y)$ can serve as a correlation kernel for the
degenerate series as well.
\endproclaim

\demo{Proof of Theorems 3.1--3.3} We prove these three theorems
simultaneously. Let $\rho_n^{(z,z',\xi)}(x_1,\dots,x_n)$ denote
the $n$--point correlation function of $M_{z,z',\xi}$. The proof
splits into two parts.

In the first part, we compute $\rho_n^{(z,z',\xi)}$ for special values of the
parameters corresponding to the degenerate series: $z=N=1,2,\dots$ and
$z'=N+\be-1$, where $\be>0$. Here we use Proposition 1.4. We show that the
formula
$$
\rho_n^{(z,z',\xi)}(x_1,\dots,x_n)
=\det[\unK_{z,z',\xi}(x_i,x_j)]_{i,j=1}^n
$$
is valid (in particular, the values of the kernel in the right--hand size are
well defined) when $z=N$, $z'=z+\be-1$, provided that $N$ is so large that the
numbers $x_i+N-\tfrac12$ are nonnegative. Then we check that in that formula,
the kernel $\unK_{z,z',\xi}$ can be replaced by the kernel $\hatK_{z,z',\xi}$:
$$
\rho_n^{(z,z',\xi)}(x_1,\dots,x_n)
=\det[\hatK_{z,z',\xi}(x_i,x_j)]_{i,j=1}^n
$$

In the second part, we extend the latter formula to other admissible values of
parameters $(z,z')$. To do this we show that both sides are analytic functions
in parameters $(z,z',\xi)$. Moreover, these functions are of such a kind that
they are uniquely determined by their values at points $(z=N, z'=N+\be-1,\xi)$.

We proceed to the detailed proof.

\proclaim{Lemma 3.4} Let $z=N=1,2,\dots$ and $z'=z+\be-1$ with $\be>0$. Assume
that $x_1,\dots,x_n$ lie in the subset $\Z_+-N+\tfrac12\subset\Z'$, so that the
points $\ti x_i:=x_i+N-\tfrac12$ are in $\Z_+$.

Then
$$
\rho_n^{(z,z',\xi)}(x_1,\dots,x_n)=\det\left[K^\m_{N,\be,\xi} (\ti x_i,\; \ti
x_j)\right]_{i,j=1}^n\,,
$$
where
$$
K^{\m}_{N,\be, \xi}(\ti x,\ti y)=\sum_{m=0}^{N-1} \wt\Me_m(\ti
x;\be,\xi)\,\wt\Me_m(\ti y;\be,\xi), \qquad  \ti x, \ti y\in\Z_+,
$$
and the functions $\wt\Me_m(\ti x;\be,\xi)$ are defined in \tht{2.8}.
\endproclaim

\demo{Proof} According to Proposition 1.3, $K^{\m}_{N,\be, \xi}(\ti x,\ti y)$
is the correlation kernel of the $N$--point Meixner orthogonal polynomial
ensemble with parameters $\be$ and $\xi$.

On the other hand, let, as above, $\Y(N)$ denote the set of Young diagrams
$\la$ with $\ell(\la)\le N$. Recall the bijective correspondence \tht{1.5}
$$
\la\mapsto \wt X(\la)=\{\ti x_1,\dots,\ti x_N\}=\{\la_1+N-1,\,
\la_2+N-2,\dots,\,\la_N\}
$$
between diagrams $\la\in\Y(N)$ and $N$--point configurations in $\Z_+$.
Comparing the definition of the infinite configuration $\unX(\la)\subset\Z'$
with that of the $N$--point configuration $\wt X(\la)$ we see that
$$
\wt X(\la)=(\unX(\la)+N-\tfrac12)\cap\Z_+.
$$

By Proposition 1.4, under this correspondence, the degenerate z--measure with
parameters $z=N$, $z'=N+\be-1$, and $\xi$ turns to the $N$--point Meixner
ensemble with parameters $\be$ and $\xi$. This implies our claim. \qed
\enddemo

We take \tht{3.3} as the definition of the kernel
$\unK_{z,z',\xi}(x,y)$.

\proclaim{Lemma 3.5} Let $z=N=1,2,\dots$ and $z'=z+\be-1$ with $\be>0$. Assume
that $x$ and $y$ lie in the subset $\Z_+-N+\tfrac12\subset\Z'$, so that $\ti
x:=x+N-\tfrac12$ and $\ti y:=y+N-\tfrac12$ are in $\Z_+$.

Then expression \tht{3.3} for the kernel $\unK_{z,z',\xi}(x,y)$
is well defined and we have
$$
\unK_{z,z',\xi}(x,y)=K^\m_{N,\be,\xi}(\ti x,\ti y).
$$
\endproclaim

\demo{Proof} We have to prove that
$$
\sum_{a\in\Z'_+}\psi_a(x;z,z',\xi)\psi_a(y;z,z',\xi) =\sum_{m=0}^{N-1}
\wt\Me_m(\ti x;\be,\xi)\,\wt\Me_m(\ti y;\be,\xi).\tag3.8
$$

We recall that the functions $\psi_a(x;z,z',\xi)$ were defined
under the assumption that both $z,z'$ are not integers. However,
as it can be seen from \tht{2.1}, each summand in the left--hand
side of \tht{3.8} makes sense under the hypotheses of the lemma.

Set
$$
a(m)=N-m-\tfrac12, \qquad m=0,1,\dots,N-1.
$$
By Proposition 2.8,
$$
\psi_{a(m)}(x;z,z',\xi)=\wt\Me_m(\ti x;\be,\xi), \qquad
\psi_{a(m)}(y;z,z',\xi)=\wt\Me_m(\ti y;\be,\xi),
$$
which implies that
$$
\sum_{a=\frac12}^{N-\frac12}\psi_a(x;z,z',\xi)\psi_a(y;z,z',\xi)
=\sum_{m=0}^{N-1} \wt\Me_m(\ti x;\be,\xi)\,\wt\Me_m(\ti y;\be,\xi). \tag3.9
$$

Finally, observe that
$$
\left.\frac1{\Ga(z-a+\tfrac12)}\,\right|_{a=N+\tfrac12,\,
N+\tfrac32,
\dots}=\left.\frac1{\Ga(N-a+\tfrac12)}\,\right|_{a=N+\tfrac12,\,
N+\tfrac32, \dots}=0.
$$

We conclude that the infinite sum in the left--hand side of
\tht{3.8} actually coincides with the finite sum in \tht{3.9}.
\qed
\enddemo

Together with Lemma 3.4 this implies

\proclaim{Corollary 3.6} Let $z=N=1,2,\dots$ and $z'=z+\be-1$ with $\be>0$.
Assume that $x_1,\dots,x_n$ lie in the subset $\Z_+-N+\tfrac12\subset\Z'$, so
that the points $\ti x_i:=x_i+N-\tfrac12$ are in $\Z_+$.

Then
$$
\rho_n^{(z,z',\xi)}(x_1,\dots,x_n)=\det\left[\unK_{z,z',\xi} (
x_i,\; x_j)\right]_{i,j=1}^n\,.
$$
\endproclaim

\proclaim{Lemma 3.7} Assume that

$\bullet$ either $(z,z')$ is not in the degenerate series  and
$x,y\in\Z'$ are arbitrary

$\bullet$ or $z=N=1,2,\dots$, $z'>N-1$, and both $x,y$ are in
$\Z_+-N+\tfrac12$.

Then the kernel $\hatK_{z,z',\xi}(x,y)$ of Theorem 3.3 is
related to the kernel $\unK_{z,z',\xi}(x,y)$ by equality
\tht{3.4}. Equivalently, the kernels are related by the ``gauge
transformation'' \tht{3.2},
$$
\hatK_{z,z',\xi}(x,y)=\frac{f_{z,z'}(x)}{f_{z,z'}(y)}\,
\unK_{z,z',\xi}(x, y), \tag3.10
$$
where $f_{z,z'}$ is defined in \tht{3.7}.
\endproclaim

\demo{Proof} Let us start with expression \tht{3.3} of the
kernel $\unK_{z,z',\xi}$ and let us replace each summand by its
integral representation \tht{2.5}. It is convenient to set
$a-\tfrac12=k$ so that as $a$ ranges over $\Z'_+$,  $k$ ranges
over $\{0,1,2,\dots\}$. Then we obtain
$$
\gather \unK_{z,z',\xi} (x,y)=\varphi_{z,z'}(x,y)\\
 \times \frac{1-\xi}{(2\pi
i)^2}\sum_{k=0}^\infty\int\limits_{\{\om_1\}}\int\limits_{\{\om_2\}}
(1-\sqrt\xi\om_1)^{-z'}
\left(1-\dfrac{\sqrt\xi}{\om_1}\right)^{z-1}
(1-\sqrt\xi\om_2)^{-z}
\left(1-\dfrac{\sqrt\xi}{\om_2}\right)^{z'-1}\\
\times \,\om_1^{-x-\tfrac12}\om_2^{-y-\tfrac12}\,
\left(\frac{(1-\sqrt\xi\om_1)(1-\sqrt\xi\om_2)}
{(\om_1-\sqrt\xi)(\om_2-\sqrt\xi)}\right)^k
\frac{d\om_1\,d\om_2}{\om_1\om_2}\,.
\endgather
$$
We can choose the contours $\{\om_1\}$ and $\{\om_2\}$ so that
they are contained in the domain $|\om|>1$. Since the
fractional--linear transformation
$$
\om\mapsto\frac{1-\sqrt\xi\om}{\om-\sqrt\xi}
$$
preserves the unit circle $|\om|=1$ and maps its exterior
$|\om|>1$  into its interior $|\om|<1$, we  have on the product
of the contours a bound of the form
$$
\left|\frac{(1-\sqrt\xi\om_1)(1-\sqrt\xi\om_2)}
{(\om_1-\sqrt\xi)(\om_2-\sqrt\xi)}\right|\le q<1.
$$
Therefore, we can interchange summation and integration and then
sum the arising geometric progression in the integrand:
$$
\sum_{k=0}^\infty
\left(\frac{(1-\sqrt\xi\om_1)(1-\sqrt\xi\om_2)}
{(\om_1-\sqrt\xi)(\om_2-\sqrt\xi)}\right)^k
=\frac{\left(1-\dfrac{\sqrt\xi}{\om_1}\right)
\left(1-\dfrac{\sqrt\xi}{\om_1}\right)\om_1\om_2}
{(1-\xi)(\om_1\om_2-1)}\,.
$$
Then we obtain equality \tht{3.4} with integral \tht{3.6}, as
desired. Finally, we can relax the assumption on the contours:
it suffices to assume that $\{\om_1^{-1}\}$ is strictly
contained inside $\{\om_2\}$, as in the formulation of Theorem
3.3.

It remains to show that \tht{3.4} is equivalent to \tht{3.10}
According to \tht{3.5} consider the expression
$$
\frac1{\varphi_{z,z'}(x,y)}=\frac{\Ga(x+z'+\tfrac12)\Ga(y+z+\tfrac12)}
{\sqrt{\Ga(x+z+\tfrac12)\Ga(x+z'+\tfrac12)
\Ga(y+z+\tfrac12)\Ga(y+z'+\tfrac12)}}
$$
Let us show that
$$
\frac1{\varphi_{z,z'}(x,y)}=\frac{f_{z,z'}(x)}{f_{z,z'}(y)}
$$
Indeed, $1/\varphi_{z,z'}$ has the form
$$
\frac{a(x)b(y)}{\sqrt{a(x)b(x)a(y)b(y)}}\,,
$$
and our hypotheses imply that $a(x)b(x)$ and $a(y)b(y)$ are real
and strictly positive. We also have
$$
f_{z,z'}(x)=\frac{a(x)}{\sqrt{a(x)b(x)}}.
$$
Therefore, we get
$$
\gather
\frac{f_{z,z'}(x)}{f_{z,z'}(y)}=\frac{a(x)\sqrt{a(y)b(y)}}{\sqrt{a(x)b(x)}\,a(y)}
=\frac{a(x)a(y)b(y)}{\sqrt{a(x)b(x)a(y)b(y)}\,a(y)}\\
=\frac{a(x)b(y)}{\sqrt{a(x)b(x)a(y)b(y)}}=\frac1{\varphi_{z,z'}(x,y)}\,.
\endgather
$$
\qed
\enddemo

\proclaim{Corollary 3.8} Let  $z=N=1,2,\dots$ and $z'>N-1$. Then
$$
\rho_n^{(z,z',\xi)}(x_1,\dots,x_n)=\det\left[\hatK_{z,z',\xi}
(x_i,\,x_j)\right]_{i,j=1}^n \tag3.11
$$
provided that all the points $x_1,\dots,x_n\in\Z'$ lie in the
subset $\Z_+-N+\tfrac12\subset\Z'$.
\endproclaim

\demo{Proof} Indeed, this follows from Lemma 3.7 and Corollary 3.6. \qed
\enddemo

This completes the first part of the proof. Now we proceed to
the second part.

\proclaim{Lemma 3.9} {\rm(i)} Fix an arbitrary set of Young diagrams $\Cal
D\subset\Y$. For any fixed admissible pair of parameters $(z,z')$, the function
$$
\xi\mapsto\sum_{\la\in \Cal D}M_{z,z',\xi}(\la),
$$
which is initially defined on the interval $(0,1)$, can be
extended to a holomorphic function in the unit disk $|\xi|<1$.

{\rm(ii)} Consider the Taylor expansion of this function at
$\xi=0$,
$$
\sum_{\la\in \Cal D}M_{z,z',\xi}(\la)=\sum_{k=0}^\infty
G_{k,\Cal D}(z,z')\xi^k.
$$
Then the coefficients $G_{k,\Cal D}(z,z')$ are polynomial
functions in $z,z'$. That is, they are restrictions of
polynomial functions to the set of admissible values $(z,z')$.
\endproclaim

\demo{Proof} (i) Set $\Cal D_n=\Cal D\cap\Y_n$. By the
definition of $M_{z,z',\xi}$,
$$
\align \sum_{\la\in \Cal D}M_{z,z',\xi}(\la)
&=\sum_{n=0}^\infty\left(\,\sum_{\la\in
\Cal D_n}M^{(n)}_{z,z'}(\la)\right) \pi_{zz',\xi}(n)\\
&=(1-\xi)^{zz'}\sum_{n=0}^\infty\left(\,\sum_{\la\in \Cal
D_n}M^{(n)}_{z,z'}(\la)\right)\frac{(zz')_n\,\xi^n}{n!}.
\endalign
$$
Each interior sum is nonnegative and does not exceed 1. On the
other hand,
$$
\sum_{n=0}^\infty |\pi_{zz',\xi}(n)|
=|1-\xi|^{zz'}\sum_{n=0}\frac{(zz')_n\,|\xi|^n}{n!}<\infty,
\qquad \xi\in\C, \quad |\xi|<1.
$$
This proves the first claim.

(ii) By \tht{1.11},
$$
\sum_{\la\in \Cal D}M_{z,z',\xi}(\la)
=(1-\xi)^{zz'}\sum_{n=0}^\infty\sum_{\la\in \Cal D_n}
(z)_\la(z')_\la\,\xi^n\left(\frac{\dim\la}{n!}\right)^2.
$$
It follows that
$$
G_{k,\Cal D}(z,z')=\sum_{n=0}^k
\frac{(-zz')_{k-n}}{(k-n)!}\sum_{\la\in \Cal D_n}
(z)_\la(z')_\la\,\left(\frac{\dim\la}{n!}\right)^2.
$$
Since each $\Cal D_n$ is a finite set, this expression is a
polynomial in $z,z'$. \qed
\enddemo

Now we can complete the proof of the theorems. Fix $n$ and an
arbitrary $n$--point subset $X=\{x_1,\dots,x_n\}\subset\Z'$, and
regard $\rho_n^{(z,z',\xi)}(x_1,\dots,x_n)$ as a function of
parameters $z,z',\xi$. We want to show that equality \tht{3.11}
holds for any admissible $(z,z')$. Apply Lemma 3.11 to the set
$\Cal D$ of those diagrams $\la$ for which $\unX(\la)$ contains
$X$, and observe that
$$
\rho_n^{(z,z',\xi)}(x_1,\dots,x_n)=\sum_{\la\in \Cal
D}M_{z,z',\xi}(\la).
$$
It follows that $\rho_n^{(z,z',\xi)}(x_1,\dots,x_n)$ is a
real--analytic function of $\xi\in(0,1)$ which admits a
holomorphic extension to the open unit disk $|\xi|<1$. Moreover,
the Taylor coefficients of this function depend on $z,z'$
polynomially.

On the other hand, from the expression \tht{3.6} for the kernel
$\hatK_{z,z',\xi}(x,y)$ it follows that this kernel (and hence
the right--hand side of \tht{3.11}) has the same property, with
$\xi$ replaced by $\sqrt\xi$.

Thus, both sides of \tht{3.11} can be viewed as (restrictions
of) holomorphic functions in $\sqrt\xi$ with polynomial Taylor
coefficients. Since the set
$$
\{(z,z')\mid \text{$z$ is a large natural number $N$ and
$z'>N-1$}\}
$$
is a set of uniqueness for polynomials in two variables, we
conclude that equality \tht{3.11} is true for any admissible
$(z,z')$.

This proves Theorem 3.1 and Theorem 3.3. Now, Theorem 3.2 follows from Theorem
3.3 and Lemma 3.7. \qed
\enddemo

\proclaim{Proposition 3.10} Formula \tht{3.3} for the kernel can also be
written as
$$
\unK_{z,z',\xi}(x,y)=\frac{\sqrt{zz'\xi}}{1-\xi}\;
\frac{\psi_{-\frac12}(x)\psi_{\frac12}(y)
-\psi_{\frac12}(x)\psi_{-\frac12}(y)}{x-y}\,. \tag3.12
$$
\endproclaim

\demo{Comment} The indeterminacy $0/0$ arising on the diagonal $x=y$ is
resolved as follows. Observe that the defining analytic expression \tht{2.1}
for $\psi_a(x)$ makes sense for any complex $x$ sufficiently close to the
lattice $\Z'$, so that we may view $\psi_{\frac12}(\,\cdot\,)$ and
$\psi_{-\frac12}(\,\cdot\,)$ as analytic functions in a neighborhood of
$\Z'\subset\C$. Since the numerator in \tht{3.13} is an analytic function in
$(x,y)$ vanishing on the diagonal $x=y$, it can be divided by $x-y$. Thus, the
value of \tht{3.12} on the diagonal can be computed, say, using the analytic
expression \tht{2.1} and the l'Hospital rule.
\enddemo

\demo{Proof} Assume first $x\ne y$. Then it suffices to prove that
$$
(x-y)\sum_{a\in\Z'_+}\psi_a(x)\psi_a(y)
=\frac{\sqrt{zz'\xi}}{1-\xi}\;(\psi_{-\frac12}(x)\psi_{\frac12}(y)
-\psi_{\frac12}(x)\psi_{-\frac12}(y)).
$$
Recall the three--term relation from Corollary 2.6, which we can write as
$$
x\psi_a(x)=A(a, a+1)\psi_{a+1}(x)+A(a,a)\psi_a(x)+A(a, a-1)\psi_{a-1}(x)
$$
with appropriate coefficients $A(\,\cdot\,,\,\cdot\,)$. Using this and the
similar relation for $y\psi_a(y)$, and taking into account the symmetry
relation
$$
A(a, a\pm1)=A(a\pm1,a)
$$
(which follows from the explicit expression in Corollary 2.6), we readily get,
after obvious cancellations,
$$
(x-y)\sum_{a\in\Z'_+}\psi_a(x)\psi_a(y)
=A(\tfrac12,-\tfrac12)(\psi_{-\frac12}(x)\psi_{\frac12}(y)
-\psi_{\frac12}(x)\psi_{-\frac12}(y)).
$$
Since
$$
A(\tfrac12,-\tfrac12)=\frac{\sqrt{zz'\xi}}{1-\xi}\,,
$$
we are done. Notice that the infinite sums involved in this computation are
convergent, because, for fixed $x$ and $y$,  $\psi_a(x)$ and $\psi_a(y)$ decay
exponentially as $a\to+\infty$: Indeed, by virtue of Proposition 2.5 this fact
reduces to that pointed out in the beginning of proof of Proposition 2.4.

To handle the case $x=y$ we use the same trick as in Lemma 3.9 and the
subsequent argument: the Taylor expansion at 0 with respect to variable
$\eta:=\sqrt\xi$.

Specifically, let us regard $\psi_a(x)=\psi_a(x;z,z',\xi)$ as a function in
$x\in\Z'$, $a\in\Z'_+$, and $\eta=\sqrt\xi$. {}From the integral representation
\tht{2.4} it is clear that this function is well defined as an analytic
function in $\eta$ ranging in the open unit disc $|\eta|<1$ in $\C$. The same
argument as above shows that this function decays exponentially as
$a\to+\infty$, uniformly on compact subsets of the disc. It follows that the
series
$$
\unK_{z,z',\xi}(x,y)=\sum_{a\in\Z'_+}\psi_a(x;z,z',\xi)\psi_a(y;z,z',\xi)
\tag3.13
$$
is analytic in the same disc $|\eta|<1$, too.

On the other hand, from \tht{2.1} it follows that if $x+a\ge0$ then
$\psi_a(x;z,z',\xi)$ is of order $O(\eta^{x+a})$ about $\eta=0$. \footnote{It
is worth noting that this claim is no longer true for negative $x+a$, because
then the hypergeometric function in the numerator of \tht{2.1} has a
singularity compensated by the gamma function in the denominator, and the order
at $\xi=0$ has to be evaluated in a more sophisticated way.} Therefore,
expanding the kernel in the Taylor series at $\eta=0$,
$$
\unK_{z,z',\xi}(x,y)=\sum_{n=0}^\infty F_n(x,y;z,z')\eta^n,
$$
we see that only finitely many terms in the series \tht{3.13} contribute to any
fixed coefficient $F_n(x,y;z,z')$. Looking again at \tht{2.1} we see that each
coefficient can be written as
$$
\multline F_n(x,y;z,z')\\
=\frac{\sqrt{(\Ga(z+x+\tfrac12)\Ga(z'+x+\tfrac12)
\Ga(z+y+\tfrac12)\Ga(z'+y+\tfrac12)}} {\Ga(x+\tfrac12)\Ga(y+\tfrac12)}\,
G_n(x,y;z,z'),
\endmultline
$$
where $G_n(x,y;z,z')$ is a {\it rational\/} function in $(x,y)$.

It follows that once we know the kernel out of the diagonal $x=y$ we can extend
it to the diagonal uniquely, by an obvious extension of the rational functions
$G_n(x,y;z,z')$. Finally, viewing the right--hand side of \tht{3.12} as an
analytic function in three variables $x$, $y$, and $\eta$, it is readily
checked that the recipe of extension suggested in the comment to the statement
of the Proposition is the correct one. \qed
\enddemo

\example{Remark 3.11} 1. The correlation functions of the z--measures
$M_{z,z',\xi}$ were first computed in \cite{BO2} in a different form: in that
paper we dealt with another embedding of partitions into the set of lattice
point configurations. The kernel $\unK_{z,z',\xi}(x,y)$ with $x,y>0$ coincides
with one of the ``blocks'' of the kernel considered in \cite{BO2}. The relation
between both kernels is discussed in detail in \cite{BO5}. The proofs in
\cite{BO2} and \cite{BO5} are very different from the arguments of the present
section.

2. Two other derivations of the kernel $\unK_{z,z',\xi}(x,y)$ are given in
Okounkov's papers \cite{Ok2} and \cite{Ok1}. In both these papers, the
correlation functions are expressed through the vacuum state expectations of
certain operators in the infinite wedge Fock space. A (substantial) difference
between the methods of \cite{Ok2} and \cite{Ok1} consists in the concrete
choice of operators. The general formalism of Schur measures presented in
\cite{Ok1} is complemented by explicit computations in \cite{BOk, \S4}. One
more derivation of the kernel $\unK_{z,z',\xi}(x,y)$ was recently suggested in
\cite{BOS}.

3. In general, kernels of the form
$$
\frac{P(x)Q(y)-Q(x)P(y)}{x-y}
$$
are called {\it integrable kernels\/}, in accordance to the terminology of
\cite{IIKS}, \cite{De1}, \cite{B}. In our case $P$ and $Q$ are expressed
through the Gauss hypergeometric function, this is why we called
$\unK_{z,z',\xi}(x,y)$ the {\it discrete hypergeometric kernel\/}.

4. The derivation of \tht{3.12} from \tht{3.3} is quite similar to the standard
derivation of the Christoffel--Darboux formula for an arbitrary system of
orthogonal polynomials. Since, as explained in \S2, the functions $\psi_a$ are
closely related to the Meixner polynomials, this similarity is not surprising.

5. Once we know that the functions $\psi_a$ form an orthonormal
basis (Proposition 2.4), the series expression \tht{3.3} for the
kernel $\unK_{z,z',\xi}(x,y)$ immediately implies that it is a
projection kernel. This fact was first proved in \cite{BO5, \S5}
in a different way.

6. The series representation \tht{3.3} is equivalent to formula
\tht{3.16} in \cite{Ok2}. A double contour representation of
various correlation kernels related to Schur measures appeared
earlier in \cite{BOk}.

7. Thus, almost all the results obtained in this section were already known.
What is really new in our paper is the approach to their derivation based on
the relationship to the Meixner polynomials. In \cite{BO7} we apply the same
approach to a more complex (dynamical) model.

8. One more novelty of the present work is appearance of the difference
operator $D$; its importance becomes especially clear in the study of the
dynamical model, see \cite{BO7}.
\endexample

\head 4. Remarks on a relationship to Krawtchouk polynomials \endhead

There exists one more possible choice of basic parameters $z$, $z'$, and $\xi$
leading to a family of probability measures on $\Y$: Namely, parameters $z,z'$
should be nonzero integers of opposite sign, while $\xi$ should be a negative
real number (thus, instead of assuming $\xi\in(0,1)$ we now require $\xi<0$).

Indeed, let $z=N$ and $z'=-N'$, where $N$ and $N'$ are positive integers, and
let $\xi<0$. Then the weight $M_{z,z',\xi}(\la)$, as defined in \tht{1.3},
vanishes unless $\ell(\la)\le N$ and $\ell(\la')\le N'$, that is, $\la$ must be
contained in the rectangle $N\times N'$. For such diagrams $\la$, we have
$(z)_\la>0$ while the sign of $(z')_\la$ equals $(-1)^{|\la|}$ (see \tht{1.2}).
Since, the sign of $\xi^{|\la|}$ also equals $(-1)^{|\la|}$, we have
$(z')_\la\xi^{|\la|}>0$. Therefore, $M_{z,z',\xi}(\la)>0$. The sum of all the
weights is still equal to 1, so that we obtain an additional family of
probability measures on $\Y$. Let us call it the {\it second degenerate
series\/} of z--measures. Its existence was pointed out in \cite{BO5, Example
1.6}.

Let $L$ be a positive integer and $p\in(0,1)$. The {\it Krawtchouk weight
function\/} with parameters $(p,L)$ is defined on the finite set
$\{0,1,\dots,L\}$  by
$$
W^\k_{p,L}(\ti x)= \binom{L}{\ti x}p^{\ti x}(1-p)^{L-\ti x}, \qquad \ti
x=0,\dots,L.
$$
The orthogonal polynomials with this weight are called the {\it Krawtchouk
polynomials\/}, see \cite{KS, \S1.10}. Let us denote them as $\frak K_n(\ti
x;p,L)$, where $n$ is the degree of the polynomial.

The next claim is a counterpart of Proposition 1.4 and can be checked directly:

\proclaim{Proposition 4.1} Under correspondence\/ \tht{1.5}, the z--measure of
the second degenerate series with parameters\/ $(z=N,\, z'=-N', \,\xi)$, where
$N,N'\in\{1,2,\dots\}$ and $\xi<0$,  turns into the $N$--point Krawtchouk
orthogonal polynomial ensemble with parameters
$$
p=\frac\xi{\xi-1}\,,\qquad L=N+N'-1.
$$
\endproclaim

Note that our assumption $\xi<0$ implies $0<p<1$.

The Krawtchouk polynomials $\frak K_n(\ti x;p,L)$ are close relatives of the
Meixner polynomials $\frak M_n(\ti x;\be,\xi)$: both families of polynomials
can be defined by the same analytic expression involving the Gauss
hypergeometric function, only the ranges of the parameters are different. The
correspondence between the families can be formally written as follows:
$$
\frak K_n(\ti x;p,L)=\frak M_n(\ti x; -L, \tfrac p{p-1}),
$$
see the very end of \S1.10 in \cite{KS}.

\proclaim{Claim 4.2} In all arguments of the present paper that rely on the
Meixner polynomials and the Meixner ensembles one could equally well use the
Krawtchouk polynomials and the corresponding ensembles.
\endproclaim

Notice that any z--measure $M_{N,-N',\xi}$ of the second degenerate series can
be written as a mixture of certain probability measures $M^{(n)}_{N,-N'}$
living on the finite sets $\Y_n$, $0\le n\le NN'$, cf. Remark 1.5. Only now the
``mixing distribution'' on $n$'s is not the negative binomial distribution but
the ordinary binomial distribution with weights
$\binom{NN'}{n}p^n(1-p)^{NN'-n}$, where $p=\frac\xi{\xi-1}$.

The measures $M^{(n)}_{N,-N'}$ admit a nice interpretation: Let
$\square_{N,N'}$ denote the rectangular Young diagram with $N$ rows and $N'$
columns. Each standard tableau $T$ of shape $\square_{N,N'}$ can be viewed as a
sequence of growing Young diagrams
$$
T=(\varnothing=\la^{(0)},\la^{(1)},\dots,\la^{(n)},\dots,
\la^{(NN')}=\square_{N,N'}),
$$
where $\la^{(n)}\in\Y_n$.

\proclaim{Proposition 4.3} Let $\{T\}$ be the set of all standard tableaux of
shape $\square_{N,N'}$ equipped with the uniform probability measure. The
push--forward of this measure under the projection $T\mapsto\la^{(n)}$
coincides with $M^{(n)}_{N,-N'}$.
\endproclaim

In this form, the measures $M^{(n)}_{N,-N'}$ appeared in \cite{PR}. A slightly
different (but essentially equivalent) interpretation can be found in
\cite{BO9, \S5}.

Notice that the measures $M^{(n)}_{z,z'}$ mentioned in Remark 1.5 can be
obtained from the measures $M^{(n)}_{N,-N'}$ by analytic continuation in the
parameters $z,z'$: this approach is developed in \cite{BO3}.

Finally, notice that the material of this section is also related to the model
considered in \cite{GTW}.

\bigskip

{\smc A.~Borodin}: Mathematics 253-37, Caltech, Pasadena, CA
91125, U.S.A.,

\medskip

E-mail address: {\tt borodin\@caltech.edu}

\bigskip

{\smc G.~Olshanski}: Dobrushin Mathematics Laboratory, Institute
for Information Transmission Problems, Bolshoy Karetny 19,
127994 Moscow GSP-4, RUSSIA.

\medskip

E-mail address: {\tt olsh\@online.ru}


\Refs

\widestnumber\key{KOV2}

\ref\key B \by A. Borodin \paper Riemann--Hilbert problem and the discrete
Bessel kernel \jour Intern. Math. Research Notices \yr 2000 \issue 9 \pages
467--494; {\tt arXiv:\,math.CO/9912093}
\endref

\ref\key BOk \by A.~Borodin and A.~Okounkov \paper  A Fredholm
determinant formula for Toeplitz determinants \jour Integral
Equations Oper. Theory \vol 37 \yr 2000 \pages 386--396; {\tt
arXiv:\, math.CA/9907165}
\endref

\ref\key BOO \by A.~Borodin, A.~Okounkov and G.~Olshanski \paper
Asymptotics of Plancherel measures for symmetric groups \jour J.
Amer. Math. Soc. \vol 13  \yr 2000 \pages 481--515; {\tt
arXiv:\, math.CO/9905032}
\endref

\ref\key BO1 \by A.~Borodin and G.~Olshanski \paper Point
processes and the infinite symmetric group \jour Math. Research
Lett. \vol 5 \yr 1998 \pages 799--816; {\tt arXiv:\,
math.RT/9810015}
\endref

\ref\key BO2 \by A.~Borodin and G.~Olshanski \paper
Distributions on partitions, point processes and the
hypergeometric kernel \jour Comm. Math. Phys. \vol 211 \yr 2000
\pages 335--358; {\tt arXiv:\, math.RT/9904010}
\endref

\ref\key BO3 \by A.~Borodin and G.~Olshanski \paper Harmonic
functions on multiplicative graphs and interpolation polynomials
\jour Electronic J. Comb. \vol 7 \yr 2000 \pages paper \#R28;
{\tt math/9912124}
\endref

\ref\key BO4 \by A.~Borodin and G.~Olshanski \paper Z--Measures
on partitions, Robinson--Schensted--Knuth correspondence, and
$\beta=2$ random matrix ensembles \inbook In: Random matrix
models and their applications (P.~Bleher and A.~Its, eds).
Cambridge University Press. Mathematical Sciences Research
Institute Publications {\bf 40}, 2001, 71--94; {\tt arXiv:\,
math.CO/9905189}
\endref

\ref\key BO5 \by A.~Borodin and G.~Olshanski \paper Random partitions and the
Gamma kernel \jour Adv. Math. \vol 194 \yr 2005 \issue 1 \pages 141--202; {\tt
arXiv:\,math-ph/0305043}
\endref

\ref\key BO6 \by A.~Borodin and G.~Olshanski \paper Z-measures on partitions
and their scaling limits \jour European Journal of Combinatorics \vol 26 \yr
2005 \issue 6 \pages 795--834; {\tt arXiv:\, math-ph/0210048}
\endref

\ref\key BO7 \by A.~Borodin and G.~Olshanski \paper Markov processes on
partitions \jour Prob. Theory and Related Fields \vol 135 \yr 2006 \issue 1
\pages 84--152
\endref

\ref\key BO8 \by A.~Borodin and G.~Olshanski \paper Stochastic dynamics related
to Plancherel measure on partitions \inbook In: Representation Theory,
Dynamical Systems, and Asymptotic Combinatorics (V.~Kaimanovich and A.~Lodkin,
eds) \publ  Amer. Math. Soc. Translations--Series 2: Advances in the
Mathematical Sciences, vol. {\bf 217}, 2006; {\tt arXiv:\, math-ph/0402064}
\endref

\ref\key BO9 \by A.~Borodin and G.~Olshanski \paper Asymptotics of
Plancherel--type random partitions \pages preprint 2006
\endref

\ref\key BOS \by A.~Borodin, G.~Olshanski, and E.~Strahov  \paper Giambelli
compatible point processes \jour Adv. Appl. Math. \vol 37 \pages 209--248; {\tt
arXiv:\, math-ph/0505021}
\endref

\ref\key De1 \by P.~Deift \paper Integrable operators \inbook In: Differential
operators and spectral theory: M. Sh. Birman's 70th anniversary collection
(V.~Buslaev, M.~Solomyak, D.~Yafaev, eds.) \bookinfo American Mathematical
Society Translations, ser. 2, v. 189 \publ Providence, R.I.: AMS \yr 1999
\pages 69--84
\endref

\ref\key De2 \by P.~Deift \book Orthogonal polynomials and random matrices: a
Riemann--Hilbert approach \bookinfo Reprint of the 1998 original \publ American
Mathematical Society \publaddr Providence, RI \yr 2000
\endref

\ref\key Er \by A.~Erdelyi (ed.) \book Higher transcendental functions. Bateman
Manuscript Project, vol. I \publ McGraw-Hill \publaddr New York \yr 1953
\endref

\ref\key GTW \by J.~Gravner, C.~A.~Tracy, and H.~Widom \paper Limit theorems
for height fluctuations in a class of discrete space and time growth models
\jour J. Statist. Phys. \vol 102  \yr 2001 \issue 5--6 \pages 1085--1132
\endref

\ref\key Gr \by F.~A.~Gr\"unbaum \paper The bispectral problem: an overview
\inbook In:  Special functions 2000: current perspective and future directions
(J.~Bustoz et al., eds). NATO Sci. Ser. II Math. Phys. Chem. \vol {\bf 30}
\pages 129--140 \publ Kluwer Acad. Publ. \publaddr Dordrecht \yr 2001
\endref

\ref\key IIKS \by A.~R.~Its, A.~G.~Izergin, V.~E.~Korepin,
N.~A.~Slavnov \paper Differential equations for quantum
correlation functions \jour Intern. J. Mod. Phys. \vol B4 \yr
1990 \pages 10037--1037
\endref

\ref\key Jo1 \by K.~Johansson \paper Shape fluctuations and random matrices
\jour Comm. Math. Phys. \vol 209 \yr 2000 \pages 437--476 {\tt arXiv:\,
math.CO/9903134}
\endref

\ref\key Jo2 \by K.~Johansson \paper Discrete orthogonal polynomial ensembles
and the Plancherel measure \jour Ann. of Math. (2) \vol 153 \yr 2001 \issue 1
\pages 259--296; {\tt arXiv:\, math.CO/9906120}
\endref

\ref\key Jo3 \by K.~Johansson \paper Non--intersecting paths, random tilings
and random matrices \jour Probab. Theory Related Fields \vol 123 \yr 2002
\issue  2 \pages 225--280; {\tt arXiv:\,math.PR/0011250}
\endref

\ref \key KOV1 \by S.~Kerov, G.~Olshanski, and A.~Vershik \paper Harmonic
analysis on the infinite symmetric group. A deformation of the regular
representation \jour Comptes Rend. Acad. Sci. Paris, S\'er. I \vol 316 \yr 1993
\pages 773--778
\endref

\ref\key KOV2 \by S.~Kerov, G.~Olshanski, and A.~Vershik \paper
Harmonic analysis on the infinite symmetric group \jour Invent.
Math. \vol 158 \yr 2004 \pages 551--642; {\tt arXiv:\,
math.RT/0312270}
\endref

\ref\key KS \by R.~Koekoek and R.~F.~Swarttouw \paper The
Askey--scheme of hypergeometric orthogonal polynomials and its
q-analogue \jour Delft University of Technology, Faculty of
Information Technology and Systems, Department of Technical
Mathematics and Informatics, Report no. 98-17, 1998 \pages
available via {\tt
http://aw.twi.tudelft.nl/$\thicksim$koekoek/askey.html}
 \endref

\ref\key K\"o \by W.~K\"onig \paper Orthogonal  polynomial ensembles in
probability theory \jour Probability Surveys \vol 2 \yr 2005 \pages 385--447
\endref

\ref\key Ma \by I.~G.~Macdonald \book Symmetric functions and Hall polynomials
\bookinfo 2nd edition \publ Oxford University Press \yr 1995
\endref

\ref\key MJD \by T.~Miwa, M.~Jimbo, E.~Date \book Solitons: Differential
equations, symmetries and infinite dimensional algebras \publ Cambridge Univ.
Press \publaddr Cambridge \yr 2000
\endref

\ref\key Ni \by S.~M.~Nishigaki \paper Level spacings at the metal--insulator
transition in the Anderson Hamiltonians and multifractal random matrix
ensembles \jour  Phys. Rev. E \vol 59 \yr 1999 \pages 2853--2862; {\tt arXiv:\,
cond-mat/9809147}
\endref

\ref\key Ok1 \by A.~Okounkov \paper Infinite wedge and measures on partitions
\jour Selecta Math. \vol 7 \yr 2001 \pages 1--25; {\tt arXiv:\,
math.RT/9907127}
\endref

\ref\key Ok2 \by A.~Okounkov \paper $SL(2)$ and $z$--measures \inbook in:
Random matrix models and their applications (P.~M.~Bleher and A.~R.~Its, eds).
Mathematical Sciences Research Institute Publications {\bf 40} \publ Cambridge
Univ. Press \yr 2001 \pages 407--420; {\tt arXiv:\,math.RT/0002136}
\endref

\ref\key Ol \by G. Olshanski \paper An introduction to harmonic analysis on the
infinite symmetric group \inbook In: Asymptotic combinatorics with applications
to mathematical physics \ed A.~M.~ Vershik \bookinfo A European mathematical
summer school held at the Euler Institute, St.~Petersburg, Russia, July 9--20,
2001 \publ Springer Lect. Notes Math. {\bf 1815}, 2003, 127--160; {\tt arXiv:\,
math.RT/0311369}
\endref

\ref\key PR \by B.~Pittel and D.~Romik \paper Limit shapes for random square
Young tableaux and plane partitions \pages {\tt arXiv:\, math.PR/0405190}
\endref

\ref\key S \by A.~Soshnikov \paper Determinantal random point fields \jour
Russian Math. Surveys  \vol 55 \yr 2000 \issue 5 \pages 923--975 (translation
from Uspekhi Mat. Nauk {\bf 55} (2000),  no. 5 (335), 107--160); {\tt arXiv:\,
math.PR/0002099}
\endref

\endRefs
\enddocument
\end